\documentclass[11pt, letterpaper]{article}
\usepackage{fontenc,amsmath,mathrsfs,amsfonts,amsthm, amssymb, MnSymbol, enumerate, hyperref, geroENG}
\usepackage{soul}
\usepackage{ulem}
\usepackage{authblk}
\usepackage{pgf,tikz}
\usepackage[top=1.25in, bottom=1.25in, left=1.0in, right=1.0in]{geometry}
\providecommand{\keywords}[1]
{
  \textbf{\textit{Keywords:}} #1
}

\providecommand{\masucla}[1]
{ 
  Mathematics Subject Classification 2020: #1	
}

\begin{document}
\date{}

\title{Hausdorff dimension of sets of numbers with large Lüroth elements}
\author{Aubin Arroyo, Gerardo González Robert}

\newcommand{\Addresses}{{
  \bigskip
  \footnotesize

  A.~Arroyo, \textsc{Instituto de Matemáticas Unidad Cuernavaca, Universidad Nacional Autónoma de México, Morelos, Mexico.}  \textit{E-mail address}: \texttt{aubinarroyo@im.unam.mx}
  \medskip
  
  G.~González Robert, \textsc{Facultad de Ciencias, Universidad Nacional Autónoma de México, Ciudad de México, Mexico.} \textit{E-mail address}: \texttt{gero@ciencias.unam.mx}
}}

\maketitle

\begin{abstract}
Lüroth series, like regular continued fractions, provide an interesting identification of real numbers with infinite sequences of integers. These sequences give deep arithmetic and measure-theoretic properties of subsets of numbers according to their growth. Although different, regular continued fractions and Lüroth series share several properties. In this paper, we explore one similarity by estimating the Hausdorff dimension of subsets of  real numbers whose Lüroth expansion grows at a definite rate. This is an extension of a result of Y. Sun and J. Wu to the context of Lüroth series. It was recently shown by Y. Feng, B. Tan, and Q.-L. Zhou that the lower bound in our main theorem is actually an equality.
\end{abstract} 


\masucla{Primary 11J83 ; Secondary 11J70, 11J99.}

\section{Introduction}

In 1883, J. Lüroth introduced a series expansion of real numbers known today as Lüroth series \cite{luroth}. He showed that for every $x\in (0,1]$ there is a unique infinite sequence $\sanu$ of natural numbers at least $2$ such that:
\begin{equation}\label{Ec01}
x = \frac{1}{a_1} + \frac{1}{a_1(a_1-1)a_2} + \ldots + \frac{1}{a_1(a_1-1)\cdots a_{n-1}(a_{n-1}-1)a_n}+ \ldots.
\end{equation}
We will write $x=\left\langle a_1,a_2,a_3,\ldots\right\rangle$ as a shorthand of \eqref{Ec01}. J. Lüroth proved that certain rational numbers might have two expansions with one of them finite and the other periodic, and that the series in \eqref{Ec01} converges for any sequence $\sanu$ in $\Na_{\geq 2}$. 

Like other expansions of real numbers, Lüroth series have been studied from a dynamical perspective (v.gr. \cite{DaKra}, \cite{jagervroedt} and the references therein). Let us describe the associated dynamical system. For any $t\in \RE$ let $[t]$ be its integer part. Following \cite{DaKra}, consider
\begin{equation}\label{DefExpDeLur}
\forall x\in[0,1]\quad
\scL(x)=
\begin{cases}
\left[ \frac{1}{x}\right]\left(\left[ \frac{1}{x} \right]+1\right)x - \left[\frac{1}{x}\right], \text{ if } 0<x\leq 1,\\
0, \quad\text{ if } x=0,
\end{cases}
\end{equation}
and consider $a_1:[0,1]\to\Na_{\geq 2}\cup\{+\infty\}$ given by $a_1(x)=[x^{-1}]+1$ for $x\in (0,1]$ and $a_1(0)=+\infty$. For $n\in\Na$ and $x\in (0,1]$ define $a_n(x)=a_1\left(\scL^n((x)\right)$, where $\scL^n:=\scL\circ\cdots\circ \scL(x)$ ($n$ copies of $\scL$). For any $x\in (0,1]$ define $(Q_n(x))_{n\geq 1}$ by
\[
Q_1(x)=a_1(x), \qquad 
\forall n\in \Na \quad Q_{n+1}(x):=Q_n(x)(a_n(x)-1)a_{n+1}(x).
\]
Let $(P_n(x))_{n\geq 1}$ be such that
\[
\forall n\in \Na\quad 
\frac{P_n(x)}{Q_n(x)}:= \frac{1}{a_1} + \frac{1}{a_1(a_1-1)a_2} + \ldots + \frac{1}{a_1(a_1-1)\cdots a_{n-1}(a_{n-1}-1)a_n};
\]
therefore, 
\[
\forall n\in \Na \quad x - \frac{P_n(x)}{Q_n(x)} = \frac{\scL^n(x)}{Q_n(x)(a_n(x)-1)},
\]
and \eqref{Ec01} holds (cfr. \cite{DaKra}, Section 1).

In analogy to regular continued fractions, we refer to $\frac{P_n(x)}{Q_n(x)}$ as the $n$\textbf{-th Lüroth convergent} and to $Q_n(x)$ as the $n$\textbf{-th Lüroth continuant}. However, unlike regular continued fractions, the positive integers $P_n(x), Q_n(x)$ might not be co-prime. We omit the dependency of $P_n$ and $Q_n$ on $x$ when there is no risk of ambiguity.

The similarities between regular continued fractions and Lüroth series go deeper. In \cite{jagervroedt}, H. Jager and C. de Vroedt showed that the Lebesgue measure on $[0,1]$, $\leb$, is $\scL$-ergodic. As a consequence, several non-trivial and well-known properties for regular continued fractions also hold for Lüroth series. For example, since $0<\leb(\{x\in [0,1]: a_1(x)=n\})<1$ for all $n\in\Na_{\geq 2}$, Birkhoff's Ergodic Theorem implies that
\[
\leb\left(\left\{ x=\left\langle a_1(x),a_2(x),\ldots,\right\rangle: \lim_{n\to\infty} a_n(x)=+\infty\right\}\right)=0.
\]
Moreover, if $\dim_H$ means Hausdorff dimension, it can be shown using elementary facts of iterated function systems (see \cite{bipe}, Theorem 2.2.2) that
\[
\dim_H \left\{ x=\left\langle a_1(x),a_2(x),\ldots,\right\rangle: \lim_{n\to\infty} a_n(x)=+\infty\right\} = \frac{1}{2}.
\]
(cfr. \cite{mun11}, Theorem 4). The previous equation is an analogue of a famous theorem by I.J. Good (\cite{good}, Theorem 1): if $x=[a_0;a_1,a_2,\ldots]$ is the regular continued fraction of a given $x\in\RE$, 
\[
\dim_H \left\{ x=[a_0;a_1,a_2,\ldots]: \lim_{n\to\infty} a_n(x)=+\infty\right\} = \frac{1}{2}.
\]
Among the vast research inspired by Good's paper \cite{good}, we can find the following  result by Y. Sun and J. Wu (\cite{suwu14}, Theorem 1.1):

\begin{teo01}[Y. Sun, J. Wu, 2014]\label{TeoSunWu}
For any $\beta>0$ we have that
\[
\dim_H \left\{ x=[0;a_1,a_2,\ldots]\in (0,1): \lim_{n\to\infty} \frac{\log a_{n+1}}{\log q_n} = \beta\right\} = \frac{1}{2+\beta}.
\]
\end{teo01}
Our main result provides a new similarity between regular continued fractions and Lüroth series, it is a weak analogue of Theorem \ref{TeoSunWu}.
\begin{teo01}[Main result]\label{TEOMAIN}
For every $\beta>0$ we have
\[
\frac{2}{3+\beta + \sqrt{\beta^2 + 6\beta +1}} \leq 
\dim_H \left\{ x=\left\langle a_1,a_2,\ldots\right\rangle : \lim_{n\to\infty} \frac{\log a_{n+1}}{\log Q_n} = \beta\right\} 
\leq \frac{1}{2 + \beta}.
\]
\end{teo01}
It was recently shown by Y. Feng, B. Tan, and Q.-L. Zhou as Theorem 1.2 in \cite{ftz} that the lower bound in Theorem \ref{TEOMAIN} is in fact an equality.

The organization of the paper is as follows. In Section 2, we introduce notation and we state some basic facts of Lüroth series. Section 3 has a proposition, Lemma \ref{LEMA01}, that will imply the lower bound of Theorem \ref{TEOMAIN}. Section 4 contains two previously known auxiliary results. We prove Theorem \ref{TEOMAIN} in Section 5. Section 6 contains the proof of Lemma \ref{LEMA01} and the Section 7 contains final remarks. By natural numbers we mean the set of positive integers and we denote them by $\Na$. We write $\Na_0:=\Na\cup\{0\}$.
\section{Notation and Elementary Facts of Lüroth Series}

From now on, $\clD:=\Na_{\geq 2}$. The map $\Lambda:\clD^\Na \to (0,1]$ associating to each sequence of integers the limit of the corresponding Lüroth series, that is $\Lambda(\sanu)=\left\langle a_1,a_2,\ldots\right\rangle$, is a continuous bijection and satisfies $\Lambda\circ\sigma=\scL\circ\Lambda$, where $\sigma$ is the left shift on $\clD^{\Na}$ (cfr. \cite{DaKra02}, Exercise 2.2.4). It is thus natural to borrow some notions from symbolic dynamics. Given $\bfa=(a_1, \ldots, a_n)\in \clD^{\Na}$, the $n$\textbf{-th level cylinder based on} $\bfa$ is
\[
\clC_n(\bfa):=\left\{x\in (0,1]: \forall j\in\{1, \ldots, n\} \quad a_j(x)=a_j\right\}.
\]
Hence, if $|\clC_n(\bfa)|$ denotes the diameter of $\clC_n(\bfa)$, we have
\[
|\clC_n(\bfa)| = \prod_{j=1}^n \frac{1}{a_j(a_j-1)}
\]
(cfr. Section 2 in \cite{jagervroedt}). For $\bfa=\sanu\in\clD^{\Na}$ and $k\in\Na$, we write $\clC_k(\bfa) = \clC_k(a_1, \ldots, a_k)$.

For each $n\in\Na$ and any $\bfa=(a_1,\ldots,a_n)\in\clD^n$, let $S_{\bfa}^n:(0,1]\to \clC_n(\bfa)$ be the inverse of $\scL^n$ restricted to $\clC_n(\bfa)$. We will sometimes write $S_{a_1,\ldots,a_n}^n$ rather than $S_{\bfa}^n$. We may extend this definition to infinite sequences in an obvious manner. For example, for $\bfa=\sanu \in \clD^{\Na}$, we have
\[
S_{\bfa}^1(x)=S_{a_1}^1(x)=\frac{1}{a_1(a_1-1)}x + \frac{1}{a_1}.
\]
It follows from the definition of $S_{\bfa}^n$ that $S_{\bfa}^n=S_{a_1}^1\circ\ldots\circ S_{a_n}^1$ and that
\[
\forall x\in (0,1]\quad 
(S_{\bfa}^{n})'(x)= \frac{1}{a_1(a_1-1)\cdots a_{n-1}(a_{n-1}-1)a_n(a_n-1)}.
\]
To make things clearer, consider $n\in\Na$, $x=\left\langle a_1,a_2,\ldots\right\rangle$, and $\bfb=(b_j)_{j\geq 1}\in\clD^{\Na}$, then
\[
\scL(x)=\left\langle a_2,a_3,a_4,\ldots \right\rangle, \quad S_{\bfb}^n(x)=\left\langle b_1,b_2,\ldots,b_n,a_1,a_2,a_3,\ldots\right\rangle.
\]

For any $n\in\Na$, $\bfa=(a_1,\ldots,a_n)\in\clD^n$ and $b\in \Na$ we write $\bfa b=(a_1,\ldots,a_n,b)$. The next proposition is a straightforward consequence of the definition of $S_{\bfa}^n$. 

In what follows, the closure of a given set $A\subseteq \RE$ is denoted by $\Cl(A)$.
\begin{propo01}\label{Propo-01}
For any $n\in\Na$, $\bfa\in \clD^n$, and $r\geq 2$ we have that
\[
\Cl \left(\bigcup_{a\geq r } \clC_{n+1}(\bfa a) \right) 
= 
\begin{cases}
S_{\bfa}^n\left([0,(r-1)^{-1}]\right), \text{ if } r\in\Na, \nonumber\\
S_{\bfa}^n\left([0,[r]^{-1}]\right), \text{ if } r\not\in\Na. \nonumber\\
\end{cases}
\]
\end{propo01}

\section{A Fundamental Lemma} 
The lower bound of Theorem \ref{TEOMAIN} will be an application of the next lemma.

\begin{lem01}\label{LEMA01}
Let $\bfs=(s_n)_{n\geq 1}$ be a sequence on $\Na_{\geq 4}$ such that $s_n\to\infty$ when $n\to\infty$. Put
\[
s_0:= \liminf_{n\to\infty} \frac{\log(s_1\cdots s_n)}{2\log(s_1\cdots s_n) + \log(s_{n+1})}. 
\]
Then,
\[
\forall N\in\Na_{\geq 2} \quad
\dim_H\left\{ x=\left\langle a_1,a_2,\ldots\right\rangle\in (0,1]: \forall n\in\Na \quad s_n\leq a_n \leq Ns_n-1\right\} = s_0.
\]
\end{lem01}

Lemma \ref{LEMA01} was originally proven for regular continued fractions in \cite{faliwawu} as Lemma 3.2. S. Munday proved an analogous result for $\alpha$-Lüroth series in \cite{mun11} (Theorem 5). While Munday's result is far more general than ours, strictly speaking, her work does not include Lemma \ref{LEMA01} (see Section \ref{SecPruebaLema} below).

\section{Auxiliary Results}
The upper bound on Theorem \ref{TEOMAIN} will be a consequence of a general estimate of the Hausdorff dimension of certain sets. Let $(X,d)$ be a complete metric space. For each non-empty set $Y\subseteq X$, $|Y|$ denotes its diameter. Let $\clA$ be a family of compact sets on $(X,d)$ with the following characteristics:
\begin{enumerate}[i.]
\item $\clA=\bigcup_{n=0}^{\infty} \clA_n$, where each $\clA_n\neq \vac$ is at most countable and $\clA_0$ contains exactly one element,
\item Every $A\in\clA$ satisfies $|A|>0$,
\item For every $n$, the members of $\clA_n$ are mutually disjoint,
\item For every $n$, for each $B \in\clA_n$ there is some $A\in\clA_{n-1}$ such that $B\subseteq A$,
\item For every $n$, for each $A \in\clA_{n-1}$ there is some $B\in\clA_{n}$ such that $B\subseteq A$,
\item $\max\{|A|:A\in\clA_n\}\to 0$ as $n\to\infty$.
\end{enumerate}

The \textbf{limit set} of $\clA$, $\mathbf{A}_{\infty}$, is 
\[
\mathbf{A}_{\infty}:= \bigcap_{n=0}^{\infty} \bigcup_{A\in\clA_n} A.
\]
In \cite{gero-good}, such a family of sets was called diametrically strongly tree-like.

The next lemma follows from the definition of $\bfA_{\infty}$ and some elementary properties of the Hausdorff dimension. A detailed proof can be found in \cite{gero-good} as Lemma 3.3. 

\begin{lem01}\label{LeGJL02}
For each $n\in\Na$ and each $A\in\clA_n$ write $D(A):=\{B\in\clA_{n+1}: B\subseteq A\}$. If $s>0$ is such that every sufficiently large $k\in\Na$ satisfies
\begin{equation}\label{EcLeGJL02}
\forall A\in\clA_k \quad
\sum_{B\in D(A)} |B|^s  \leq |A|^s,
\end{equation}
then $\dim_H \mathbf{A}_{\infty}\leq s$.
\end{lem01}

We recall the basic yet important Mass Distribution Principle (see \cite{bipe}, Lemma 1.2.8 for a proof).

\begin{lem01}[Mass Distribution Principle]
Let $(X,d)$ be a metric space and for each $x\in X$ and $r>0$ let $B(x;r)$ be the open ball of center $x$ and radius $r$. If $E\subset X$ supports a strictly positive Borel measure, $\mu$, such that there exists $C>0$ and $\delta>0$ verifying
\[
\forall x\in X \quad \forall r>0 \quad \mu(B(x;r))\leq C r^{\delta},
\]
then $\dim_H E\geq \delta$.
\end{lem01}
\section{Proof of Main Theorem}
Define the family of sets $\{G(\alpha): \alpha>0\}$ by
\[
\forall \alpha>0\quad 
G(\alpha) := \left\{ x \in (0,1] : \lim_{n\to\infty} \frac{\log Q_n(x)}{\log a_{n+1}(x)} = \alpha\right\}.
\]
Hence, Theorem \ref{TEOMAIN} is equivalent to
\begin{equation}\label{Eq-MainResult}
\forall \alpha>0 \quad
\frac{2\alpha}{1+3\alpha+\sqrt{\alpha^2+6\alpha+1}} \leq \dim_H G(\alpha) \leq \frac{\alpha}{2\alpha +1}.
\end{equation}
\subsection{Lower Bound in Theorem \ref{TEOMAIN}}
In this section, we study the Hausdorff dimension of sets of numbers whose Lüroth sequence $\sanu$ grows in a controlled way. Afterwards, we use this computation to obtain the lower bound in Theorem \ref{TEOMAIN}. In what follows, $c^{b^a}:=c^{(b^a)}$ for $a,b$, and $c>0$.
\begin{def01}
For any $c,d\in\clD$ and any $\lambda>1$ define
\[
K_d^c(\lambda) = \left\{ x \in (0,1]: \forall n\in\Na \quad c^{\lambda^n} \leq a_n(x) < d\, c^{ \lambda^n }\right\}.
\]
\end{def01}

\begin{lem01}\label{Le-LoBo01}
Let $\lambda>1$ and $c,d\in\clD$ be arbitrary. Every $x=\left\langle a_1,a_2,\ldots\right\rangle\in K_d^c(\lambda)$ satisfies
\begin{equation}\label{Eq-LoBo01}
\lim_{n\to\infty} \frac{\log Q_n}{\log a_{n+1}} = \frac{\lambda + 1}{\lambda(\lambda-1)}.
\end{equation}
\end{lem01}
\begin{proof}
Keep the statement's notation. Let $n$ be any natural number. Recalling the definition of $Q_n$, we have that
\begin{align*}
\log Q_n &= \log\left( a_1^2\ldots a_{n-1}^2 a_n \left(1 - \frac{1}{a_{1}}\right)\cdots \left(1 - \frac{1}{a_{n-1}}\right) \right)\nonumber\\
&\leq \log a_{n} + 2\sum_{j=1}^{n-1} \log a_j \nonumber\\
&\leq \lambda^n \log c + \log d + 2\sum_{j=1}^{n-1} (\lambda^j \log c + \log d) \nonumber\\
&= \lambda^n \log c \left( 1+ \frac{\log d}{\lambda^n \log c} + \frac{2(n-1)\log d}{\lambda^n \log c} + \frac{2}{\lambda-1} \left( 1-\frac{1}{\lambda^{n-1}} \right)\right),\nonumber
\end{align*}
and, since $\log a_{n+1}\geq  \lambda^{n+1}\log  c$, 
\begin{equation}\label{Eq-LoBo02}
\frac{\log Q_n}{\log a_{n+1}} \leq \frac{1}{\lambda}\left( 1+ \frac{\log d}{\lambda^n \log c} + \frac{2(n-1)\log d}{\lambda^n \log c} + \frac{2}{\lambda-1} \left( 1-\frac{1}{\lambda^{n-1}} \right)\right).
\end{equation}
On the other hand, the definition of $Q_n$ and $\lambda^n\log c \leq \log a_n$ give
\begin{align*}
\log Q_n &\geq \log a_n + 2\sum_{j=1}^{n-1}\log a_j -(n-1)\log 2 \nonumber\\
&\geq \lambda^n \log c + 2\sum_{j=1}^{n-1} \lambda^j \log c - (n-1)\log 2 \nonumber\\
&= \lambda^n \log c \left( 1+ \frac{2}{\lambda-1}\left( 1- \frac{1}{\lambda^{n-1}}\right) - \frac{(n-1)\log 2}{\lambda^n \log c}\right).\nonumber
\end{align*}
In view of $\log a_{n+1} \leq \lambda^{n+1} \log c + \log d$, the previous inequalities imply
\begin{equation}\label{Eq-LoBo03}
\frac{\log Q_n}{\log a_{n+1}}\geq \frac{1}{\lambda + \frac{\log d}{\lambda^n \log c}} \left( 1+ \frac{2}{\lambda-1} \left( 1- \frac{1}{\lambda^{n-1}}\right) - \frac{(n-1)\log 2}{\lambda^n \log c}\right).
\end{equation}
Letting $n\to\infty$ in \eqref{Eq-LoBo02} and \eqref{Eq-LoBo03} we conclude \eqref{Eq-LoBo01}.
\end{proof}
\begin{lem01}\label{Le-LoBo02}
For every $c\in \Na_{\geq 4}$, $d\in\clD$ and any $\lambda>1$ we have
\[
\dim_H \, K_d^c (\lambda) = \frac{1}{1+\lambda}.
\]
\end{lem01}

\begin{proof}
Let $c,d,\lambda$ be as in the statement and define $(s_n)_{n\geq 1}$ by $s_n=c^{\lambda^n}$. Then, for $n\in\Na$
\[
\log(s_1\cdots s_n) = \log c  \sum_{j=1}^n \lambda^j = \frac{\lambda^{n+1}\log c}{\lambda-1} \left(1-\frac{1}{\lambda^n}\right),
\]
and
\begin{align*}
2\log(s_1\cdots s_n) + \log(s_{n+1} ) &= \frac{\lambda^{n+1}\log c}{\lambda-1} \left(2-\frac{2}{\lambda^n}\right) + \lambda^{n+1} \log c  \nonumber\\
&= \frac{\lambda^{n+1}\log c}{\lambda-1} \left( 1 + \lambda - \frac{2}{\lambda^n} \right). \nonumber
\end{align*}
As a consequence, 
\[
\frac{\log(s_1\cdots s_n)}{2\log(s_1\cdots s_n) + \log(s_{n+1} )} = 
\frac{1 - \lambda^{-n}}{1+\lambda - 2\lambda^{-n}} \to \frac{1}{1+\lambda} \quad\text{ as }\quad n\to\infty, 
\]
and, by Lemma \ref{LEMA01}, $\dim_H K_d^c(\lambda)= (1+\lambda)^{-1}$.
\end{proof}

\begin{proof}[Proof of the lower bound in Theorem \ref{TEOMAIN}]
Take $\alpha>0$ and let $\lambda>1$ be such that $\alpha= \tfrac{\lambda+1}{\lambda(\lambda-1)}$; that is,
\[
\lambda = \frac{1+\alpha + \sqrt{\alpha^2 + 6\alpha+1}}{2\alpha}>1.
\]
Because of Lemma \ref{Le-LoBo01}, such $\lambda$ and any $c,d\in\Na_{\geq 2}$ verify $K_{c}^d(\lambda)\subseteq G(\alpha)$. Then, Lemma \ref{Le-LoBo02} implies
\begin{equation}\label{Eq-LowerBound}
\dim_H G(\alpha) \geq \dim_H K_d^c(\lambda) = \frac{1}{1+\lambda} = \frac{2\alpha}{1+3\alpha+\sqrt{\alpha^2+6\alpha+1}}. \quad \qedhere
\end{equation}
\end{proof}
We shall require the following technical result. Although it was already used in \cite{faliwawu}, we add its proof for completeness' sake.
\begin{lem01}
Let $(s_n)_{n\geq 1}$ be a sequence of natural numbers such that $s_n\to\infty$ when $n\to\infty$, then
\begin{equation}\label{Ec-s-04}
\lim_{n\to\infty} \frac{\log(s_1\cdots s_n)}{n} = + \infty.
\end{equation}
Furthermore, if we define
\begin{equation}\label{Ec-s-05}
s_0:=\liminf_{n\to\infty} \frac{\log(s_1\cdots s_n)}{2\log(s_1\cdots s_n) + \log s_{n+1}},
\end{equation}
then, for every $0<s<s_0$ and $N\in\Na_{\geq 2}$ there exists $n_0\in\Na$ such that
\begin{equation}\label{Ec-s-06}
\forall n\in\Na_{\geq n_0} \qquad \left( s_{n+1} \left( \prod_{k=1}^n Ns_k\right)^2\right)^s \leq \prod_{k=1}^n (N-1)s_k.
\end{equation}
\end{lem01}
\begin{proof}
Keep the statement's notation. Let $c>0$ be arbitrary and let $N\in\Na$ be such that $s_{N+n}\geq e^c$ for any $n\in\Na_0$. Then, for any natural number $n$ we have $s_1\cdots s_Ns_{N+1}\cdots s_{N+n}\geq e^{c(n+1)}$ and thus
\[
\frac{\log(s_1\cdots s_{N+n})}{N+n} = \frac{\log(s_1\cdots s_{N+n})}{n+1} \,\frac{n+1}{n+N}\geq c \frac{n+1}{n+N}.
\]
We conclude \eqref{Ec-s-04} by taking the inferior limit over $n$.

In order to prove the second part, let $n$ be any natural number. Then, $n$ satisfies \eqref{Ec-s-06} if and only if
\begin{equation}\label{EcLem42Aux}
s\leq  \frac{n\log(N-1) + \log(s_1\cdots s_n)}{2n\log(N) + 2 \log(s_1\cdots s_n) + \log(s_{n+1})}.
\end{equation}
Define 
\[
\alpha_n = \frac{\log(N-1)}{2\log(N) + 2\frac{\log(s_1\cdots s_n)}{n} + \frac{\log(s_{n+1})}{n}}, \qquad
\beta_n = \frac{\log(N)}{2\log(N) + 2\frac{\log(s_1\cdots s_n)}{n} + \frac{\log(s_{n+1})}{n}}.
\]
Then, we can rewrite the right-hand side in \eqref{EcLem42Aux} as
\[
\rho_n:=\alpha_n - 2\left(\frac{\log(s_1\cdots s_n)}{2\log(s_1\cdots s_n) + \log(s_{n+1})}\right)\beta_n + \frac{\log(s_1\cdots s_n)}{2\log(s_1\cdots s_n) + \log(s_{n+1})}.
\]
The coefficient of $\beta_n$ in the above expresion is bounded and, by \eqref{Ec-s-04}, $\alpha_n\to 0$ and $\beta_n\to 0$ as $n\to\infty$. Then, $\liminf_n \rho_n=s_0>s$ and \eqref{EcLem42Aux} holds whenever $n$ is large enough.
\end{proof}

\subsection{Upper Bound in Theorem \ref{TEOMAIN}}

This section is devoted to prove the upper bound of Theorem \ref{TEOMAIN}. First, we  introduce some notation. Then, we prove an auxiliary result (Proposition \ref{Prop-UsGJL}) and afterwards we use it to conclude the upper bound in \eqref{Eq-MainResult}.

Fix $\alpha>0$. Take $\veps >0$ and define: 
\begin{align*}
C_{\veps}(\alpha) &:= \left\{ x \in (0,1] : \limsup_{n\to\infty} \frac{\log Q_n(x)}{\log a_{n+1}(x)} \leq \alpha + \frac{\veps}{2}\right\}, \nonumber\\
\forall n\in\Na \qquad C_{\veps}^n(\alpha) &:= \left\{ x \in (0,1] : \forall m\in\Na_{\geq n} \; \frac{\log Q_{m}(x)}{\log a_{m+1}(x)} < \alpha + \veps \right\}.
\end{align*}
Therefore, in view of $G(\alpha)\subseteq C_{\veps}(\alpha) \subseteq \bigcup_{n\in\Na} C_{\veps}^n(\alpha)$, we have
\begin{equation}\label{EqHausDimGalfa}
\dim_H G(\alpha)
\leq \dim_H C_{\veps}(\alpha) 
\leq \sup_{n\in\Na} \; \dim_H C_{\veps}^n(\alpha) 
= \lim_{n\to\infty} \dim_H C_{\veps}^n(\alpha).
\end{equation} 
Aiming towards a lighter notation, for any $n\in\Na$ and any $\bfb$ in $\clD^{n}$ we write 
\[
R_n:=R_n(\bfb;\alpha,\veps):=Q_n(\bfb)^{\frac{1}{\alpha+\veps}}. 
\]
For each $\bfb\in\clD^n$ consider
\begin{equation}\label{Eq-DimBo}
K_n(\bfb):= K_n(\bfb;\alpha,\veps) = \Cl   \left(\bigcup \left\{\clC(\bfb b_{n+1}) : b_{n+1} \geq R_n(\bfb) \right\} \right) .
\end{equation}

By Proposition \ref{Propo-01}, $K_n(\bfb)$ is a compact interval and satisfies the following bounds
\begin{equation}\label{Eq-BoDiam}
\frac{1}{R_n(\bfb)}\prod_{j=1}^n \frac{1}{b_j(b_j-1)} 
\leq
|K_n(\bfb)|
\leq 
\frac{1}{R_n(\bfb)-1}
\prod_{j=1}^n \frac{1}{b_j(b_j-1)}.
\end{equation}
\begin{propo01}\label{Prop-UsGJL}
For each $\delta_1$ satisfying $0<\delta_1< (2+1/(\alpha+\veps))^{-1}$ let $s=s(\veps,\delta_1)= \frac{\alpha+ \veps}{2(\alpha + \veps)+1} + \delta_1$. There exists $N=N(\alpha,\veps,\delta_1)\in\Na$ such that
\[
\forall n\in\Na_{\geq N} \quad \forall \bfb\in \clD^n \quad
\sum_{c\geq R_n(\bfb)} |K_{n+1}(\bfb c)|^s \leq |K_n(\bfb)|^s.
\]
\end{propo01}
\begin{proof}
Let $\delta_1,s$ be as in the statement and define $\delta=\delta_1\left( 2 + \frac{1}{\alpha + \veps}\right)\in (0,1)$. Let $N\in\Na$ be such that
\[
N\geq \max\left\{\alpha + \veps, \frac{(\alpha+\veps)}{\delta} \log_2 \left(\frac{8}{\delta}\right)\right\}.
\]
Take $n\in\Na_{\geq N}$, $\bfb\in \clD^n$, and $c\in \Na$, $c\geq R_{n}(\bfb)$. Since $Q_m(x)\geq 2^m$ for every $x\in (0,1]$ and $m\in\Na$, the choice of $N\geq \alpha + \veps$ implies $R_{n+1}(\bfb c)>R_{n}(\bfb)>2$, so $1-c^{-1}>2^{-1}$ and $1-R_{n+1}(\bfb c)^{-1}>2^{-1}$. Direct computations and $b_n\geq 2$ give
\[
R_{n+1}(\bfb c)-1
=
c^{\frac{1}{\alpha + \veps}} (b_n - 1 )^{\frac{1}{\alpha + \veps}} R_n(\bfb)\left( 1- R_{n+1}(\bfb c)^{-1}\right) > \frac{1}{2} c^{\frac{1}{\alpha+\veps}} R_n(\bfb).
\]
Thus, the bounds in \eqref{Eq-BoDiam} imply
\begin{align*}
|K_{n+1}(\bfb c)| &\leq \frac{1}{(R_{n+1}(\bfb c)-1)c(c-1)} \prod_{j=1}^n \frac{1}{b_j(b_j-1)} \nonumber\\
&= \frac{2}{c^{ 2+\frac{1}{\alpha + \veps} } R_{n}(\bfb)(1-c^{-1})} \prod_{j=1}^n \frac{1}{b_j(b_j-1)} \nonumber\\
&\leq \frac{4}{c^{ 2+\frac{1}{\alpha + \veps} }}\, \frac{1}{R_{n}(\bfb)} \prod_{j=1}^n \frac{1}{b_j(b_j-1)} \nonumber\\
&\leq 4 c^{-\left( 2+\frac{1}{\alpha + \veps} \right)} |K_n(\bfb)|. 
\end{align*}
Varying $c$ along the integers larger than or equal to $R_n=R_n(\bfb)$ we obtain
\begin{align*}
\sum_{c \geq R_n} |K_{n+1}(\bfb c)|^s &\leq 4^s |K_n(\bfb)|^s \sum_{c \geq R_n} c^{-(1+\delta)} \nonumber\\
&\leq 4 |K_n(\bfb)|^s \int_{R_n-1}^{\infty} x^{-(1+\delta)}\md x \nonumber\\
&= \frac{4}{\delta R_n^{\delta}(1-R_n^{-1})^{\delta}} |K_n(\bfb)|^s\nonumber\\
&< \frac{4\cdot 2^{\delta}}{\delta R_n^{\delta}} |K_n(\bfb)|^s\nonumber\\
&\leq \frac{8}{\delta R_n^{\delta}} |K_n(\bfb)|^s \leq |K_n(\bfb)|^s. \nonumber
\end{align*}
The last inequality follows from $n\geq N\geq \tfrac{\alpha+\veps}{\delta}\log_2 \left(\tfrac{8}{\delta}\right)$.

\end{proof}

Now we can obtain the desired upper bound in \eqref{Eq-MainResult}. Define the following subsets of natural numbers
\[
\forall n\in\Na \quad \forall \bfa\in \clD^n\quad
I(\bfa):=I(\bfa,n;\alpha,\veps) := \{c\in\Na: c \geq R_n(\bfa)\}.
\]
Consider a positive $\delta_1<(2+1/(\alpha + \veps))^{-1}$ and let $s,N$ be as in Proposition \ref{Prop-UsGJL}. Take $\bfb\in\clD^N$. Define $\{\clA_k:k\in\Na_0\}$ by $\clA_0  := \{ K_N(\bfb)\}$ and
\[
\forall j\in\Na \quad 
\clA_j :=
\{K_{N+j}(\bfb c_1\cdots c_j): c_1\in I(\bfb), c_2\in I(\bfb c_1), \ldots, c_{n+j}\in I(\bfb c_1\cdots c_{j-1})\}. \nonumber
\]
We can apply Lemma \ref{LeGJL02} on $\{\clA_k:k\in\Na_0\}$ to conclude that the Hausdorff dimension of its limit set is at most $s$:
\[
\dim_H \; \bigcap_{n=0}^{\infty} \bigcup_{A\in\clA_n} A\; =  \dim_H \; C_{\veps}^N(\alpha)\cap K_N(\bfb)\; \leq s.
\]
Since $\bfb$ was arbitrary and
\[
C_{\veps}^N(\alpha) = \bigcup  \left\{ C_{\veps}^N(\alpha)\cap K_N(\bfa): \bfa\in \clD^N\right\},
\]
we obtain $\dim_H C_{\veps}^N(\alpha)\leq s$. Furthermore, since $\delta_1>0$ can be arbitrarily small, from \eqref{EqHausDimGalfa} we obtain
\[
\dim_H G(\alpha) \leq \frac{\alpha + \veps}{2(\alpha + \veps) + 1}.
\]
Finally, the function $x\mapsto \frac{x}{2x +1}$ is continuous and strictly increasing on $\RE_{>0}$, so we get $\dim_H G(\alpha) \leq \frac{\alpha}{2\alpha+1}$ by letting $\veps\to 0$. This concludes the proof of the upper bound in \eqref{Eq-MainResult}.

\section{Proof of Lemma \ref{LEMA01}}\label{SecPruebaLema}

Before showing Lemma \ref{LEMA01}, we elaborate on our earlier claim that S. Munday's work in \cite{mun11} does not include it. Let us recall the definition of $\alpha$-Lüroth series. Let $\alpha=\{I_n:n\in\Na\}$ be a partition of $[0,1)$ where each $I_n$ is left closed and right open, except for $I_1$, which is open. Suppose that
\[
\forall m,n\in\Na \quad m<n \implies \sup I_n \leq \inf I_m. 
\]
For each $n\in\Na$, write $i_n=|I_n|$ and $t_n=\sum_{k\geq n} i_n$. Define
\[
\forall x\in [0,1) \quad 
\scL_{\alpha}(x) = 
\begin{cases}
(t_n-x) i_n^{-1}, \;\text{if}\; x\in I_n,\\
0, \;\text{if}\; x=0.
\end{cases}
\]
Using $\scL_{\alpha}$, we can associate to each $x\in [0,1)$ a sequence $\bfb=\sabu$ on $\Na$ using the condition $\scL_{\alpha}^{n-1}(x)\in I_{b_j}$ for every $j$. Such sequence satisfies
\[
x= t_{b_1} - i_{b_1}t_{b_2} + i_{b_1}i_{b_2}t_{b_3} - i_{b_1}i_{b_2}i_{b_3}t_{b_4} + \ldots.
\]
This expansion is the \textbf{$\alpha$-Lüroth} series of $x$. Under additional conditions, S. Munday showed in \cite{mun11} that Lemma \ref{LEMA01} holds for $\alpha$-Lüroth series.  While $\alpha$-Lüroth series are alternating, Lüroth series are not. Although it might sound irrelevant, this slight modification gives rise to new complications when working out the details (see the third case in the proof of Proposition \ref{PropoUniqueInterFI}).


\subsection{Terminology and notation}

Let $\bfs=(s_n)_{n\geq 1}$ be a sequence of natural numbers such that $s_n\geq 4$ for all $n\in\Na$, and $s_n\to +\infty$ when $n\to \infty$ and take $N\in\Na_{\geq 2}$. Define
\begin{align*}
\forall n\in\Na \quad \clJ(\bfs,n) :=\left\{ \bfa=(a_1,\ldots, a_n)\in\Na^n: \forall j\in\{1,\ldots,n\} \quad s_j\leq a_j \leq N s_j -1 \right\},\nonumber\\
\forall n\in\Na \quad \forall \bfa=(a_1,\ldots,a_n) \in \clJ(\bfs,n) \quad
J_n(\bfa):= \Cl\left(\bigcup_{a_{n+1}\geq s_{n+1}} \clC_{n+1}(\bfa a_{n+1})\right).
\end{align*}

Hence, every $n\in\Na$ satisfies
\begin{equation}\label{DiamJna}
\forall \bfa=(a_1,\ldots,a_n)\in\clJ(\bfs,n) 
\qquad 
|J_n(\bfa)|= \frac{1}{s_{n+1}-1}\prod_{j=1}^n \frac{1}{a_j(a_j-1)}.
\end{equation}
As we did for cylinders, we extend the notation $J_n(\bfa)$ to infinite sequences $\bfa$. In this case, we have $\clC_{n+1}(\bfa)\subseteq J_n(\bfa)\subseteq \clC_n(\bfa)$ for all $n\in\Na$. We call the sets $J_n(\bfa)$ \textbf{fundamental intervals} and by the \textbf{order} of $J_n(\bfa)$ we mean $n$. A direct consequence of \eqref{DiamJna} is that every $n\in\Na$ and $\bfa\in \clJ(\bfs,n)$ satisfy
\begin{equation}\label{Ec-CotaDiamJna}
\frac{1}{N^{2n+1}(s_1\cdots s_n)^2s_{n+1}} 
\leq |J_n(\bfa)| \leq
\frac{2^n}{(s_1\cdots s_n)^2 s_{n+1}}.
\end{equation}
For the rest of the proof, we will denote by $F$ the set in the statement of Lemma \ref{LEMA01}:
\[
F:=F(\bfs):=\left\{ x=(0,1]: \forall n\in \Na \quad s_n\leq a_n(x) \leq Ns_n-1\right\}.
\]
Our goal, then, is to compute the Hausdorff dimension of $F$.

\subsection{Upper Bound of $\dim_H F$}

We cover $F$ with fundamental intervals to conclude an upper estimate of $\dim_H F$. Define
\[
s_0:= \liminf_{n\to\infty} \frac{\log(s_1\cdots s_n)}{2\log(s_1\cdots s_n) + \log(s_{n+1})} 
\]
and take $s$ with $s_0< s\leq 1$. By definition of $s_0$, there is an increasing sequence of natural numbers $(n_j)_{j\geq 1}$ such that
\begin{equation}\label{EcCotaSup01}
\forall j\in\Na \quad \prod_{k=1}^{n_j} s_j\leq \left(s_{n_j+1}(s_1\cdots s_{n_j})^2\right)^{\tfrac{s+s_0}{2}}.
\end{equation}
In view of \eqref{Ec-s-04}, we may assume that $(n_j)_{j\geq 1}$ also satisfies
\[
\forall j\in \Na \quad 2^{n_j}(N-1)^{n_j}< (s_{n_j+1}(s_1\cdots s_{n_j})^2)^{\tfrac{s-s_0}{2}}.
\]
Hence, for all $j\in\Na$
\begin{align*}
\sum_{\bfa\in \clJ(\bfs,n_j)} |J_{n_j}(\bfa)|^s &\leq (N-1)^{n_j} \prod_{k=1}^{n_j} s_k \left( \frac{2^{n_j}}{(s_1\cdots s_{n_j})^2 s_{n_j+1}} \right)^s \nonumber\\
&< \frac{ \left((s_1 \cdots s_{n_j})^2 s_{n_j+1} \right)^{\frac{s-s_0}{2}}}{2^{n_j}}\left( \prod_{k=1}^{n_j} s_k\right) \frac{2^{sn_j}}{\left( (s_1\cdots s_{n_j})^{2} s_{n_j+1}\right)^s} \nonumber\\
&= \left(\frac{1}{2^{1-s}}\right)^{n_j} \frac{1}{\left((s_1\cdots s_{n_j})^2s_{n_j+1} \right)^{\frac{s+s_0}{2}}} \prod_{k=1}^{n_j} s_k \leq \left(\frac{1}{2^{1-s}}\right)^{n_j} \leq 1.
\end{align*}

Finally, we estimate the $s$-Hausdorff measure of $F$, $\clH^s(F)$, as follows:
\[
\clH^s(F)\leq \liminf_{n\to\infty} \sum_{\bfa\in \clJ(\bfs,n)} |J_n(\bfa)|^s \leq \liminf_{j\to\infty} \sum_{\bfa\in \clJ(\bfs,n_j)} |J_{n_j}(\bfa)|^s \leq 1,
\]
so $\dim_H F \leq s$. Since $s\geq s_0$ was arbitrary, $\dim_H F \leq s_0$.

\subsection{Lower Bound of $\dim_H F $}
We start by defining a measure $\tilde{\mu}$ on the space $\prod_{n\geq 1} \{s_n,\ldots, Ns_n-1\}$ equipped with the $\sigma$-algebra induced by the product topology assuming each factor has the discrete topology. For each $n\in\Na$, define the measure $\tilde{\mu}_n$ on $\{s_n,\ldots, Ns_n-1\}$ by
\[
\forall k\in \{s_n,\ldots, Ns_n-1\} \quad \tilde{\mu}_n \left(\{k\}\right)=\frac{1}{(N-1)s_n}.
\]
Using the Daniell-Kolmogorov Consistency Theorem (\cite{parth}, Proposition 3.6.4), we obtain the measure $\tilde{\mu}$ for the sequence $(\{s_n,\ldots,Ns_n-1\},\tilde{\mu}_n)_{n\geq 1}$. Let $\Lambda:\clD^{\Na} \to [0,1)$ be given by $\sanu\mapsto \left\langle a_1,a_2,\ldots\right\rangle$ and define $\mu=\tilde{\mu} \Lambda^{-1}$. Hence,
\[
\forall n\in\Na \quad \forall \bfa\in\clJ(\bfs,n) \quad \mu(J_n(\bfa))=\prod_{j=1}^n \frac{1}{(N-1)s_j}.
\]
Take $0<s<s_0$, let $n_0\in\Na$ be as in \eqref{Ec-s-06} and $r_0:=(N^{2{n_0}+1}(s_1\cdots s_{n_0})^2s_{{n_0}+1})^{-1}$. 

We proceed to verify that $\mu$ satisfies the hypotheses of the Mass Distribution Principle. Take $x=\left\langle a_1,a_2,\ldots\right\rangle\in F$ and $0<r<r_0$. Write $\bfa=\sanu$ and call $n$ the unique natural number larger than or equal to $n_0$ such that
\[
|J_{n+1}(\bfa)|\leq r < |J_n(\bfa)|.
\]
The following proposition, whose proof we delay to the end of this section, will be heavily used.
\begin{propo01}\label{PropoUniqueInterFI}
The only fundamental interval of order $n$ intersecting $B(x;r)$ is $J_n(\bfa)$.
\end{propo01}

We consider two cases to prove the existence of some $C=C(N)>0$ such that $\mu(B(x;r))\leq C r^{s}$ holds. These cases are:
\[
|J_{n+1}(\bfa)|\leq r <|\clC_{n+1}(\bfa)|, \quad  |\clC_{n+1}(\bfa)|\leq r<|J_n(\bfa)|.
\]

\paragraph{Case I.} Assume that $|J_{n+1}(\bfa)|\leq r <|\clC_{n+1}(\bfa)|$. We claim that
\begin{equation}\label{Ec-Propo53}
B(x;|\clC_{n+1}(a_1,\ldots,a_{n+1})|) \subseteq \bigcup_{j=-1}^2 \clC_{n+1}(a_1,\ldots,a_n, a_{n+1}+j).
\end{equation}
Indeed, on the one hand, from
\begin{align*}
|\clC_{n+1}(a_1,\ldots,a_n, a_{n+1}-1)| &= \left( \prod_{j=1}^{n} \frac{1}{a_j(a_j-1)}\right) \frac{1}{(a_{n+1}-1)(a_{n+1}-2)} \nonumber\\
&=  \frac{a_{n+1}}{a_{n+1}-2}\,\left|\clC_{n+1}(a_1,\ldots,a_{n+1})\right|>|\clC_{n+1}(a_1,\ldots,a_{n+1})| \nonumber
\end{align*}
and $\sup\clC_{n+1}(a_1,\ldots,a_{n+1})=\inf\clC_{n+1}(a_1,\ldots,a_n, a_{n+1}-1)$, we have
\begin{align}
x+|\clC_{n+1}(a_1,\ldots,a_{n+1})| &\leq \sup \clC_{n+1}(a_1,\ldots,a_{n+1}) + |\clC_{n+1}(a_1,\ldots,a_{n+1})| \nonumber\\
&\leq \inf \clC_{n+1}(a_1,\ldots,a_n, a_{n+1}-1) + |\clC_{n+1}(a_1,\ldots,a_n,a_{n+1}-1)| \nonumber\\
&=\sup \clC_{n+1}(a_1,\ldots,a_n, a_{n+1}-1). \label{Ec-PropoAux-01}
\end{align}
On the other hand, we have that
\begin{align*}
\inf\clC_{n+1}(a_1,\ldots,a_{n+1}) - \inf\clC_{n+1}(a_1,\ldots,a_n, a_{n+1}+2)  
&= S_{a_1}\cdots S_{a_{n}}(S_{a_{n+1}}(0)) - S_{a_1}\cdots S_{a_{n}}(S_{a_{n+1}+2}(0)) \nonumber\\
&= \left( \prod_{j=1}^{n} \frac{1}{a_j(a_j-1)}\right) \left( \frac{1}{a_{n+1}} - \frac{1}{a_{n+1}+2}\right) \nonumber\\
&= \left( \prod_{j=1}^{n} \frac{1}{a_j(a_j-1)}\right) \left( \frac{2}{a_{n+1}(a_{n+1}+2)} \right) \nonumber\\
&\geq |\clC_{n+1}(a_1,\ldots,a_n, a_{n+1})|.
\end{align*}
The last equality follows from $a_{n+1}\geq 4$. Hence, 
\[
\inf\clC_{n+1}(a_1,\ldots,a_{n+1})-|\clC_{n+1}(a_1,\ldots,a_{n+1})|\geq \inf\clC_{n+1}(a_1,\ldots,a_n, a_{n+1}+2) 
\]
and, since $x>\inf \clC_{n+1}(a_1,\ldots, a_{n+1})$, we have that
\begin{equation}\label{Ec-PropoAux-02}
x - |\clC_{n+1}(a_1,\ldots,a_{n+1})| > \inf\clC_{n+1}(a_1,\ldots,a_n, a_{n+1}+2).
\end{equation}
The inequalities \eqref{Ec-PropoAux-01} and \eqref{Ec-PropoAux-02} imply \eqref{Ec-Propo53}.

With Proposition \ref{PropoUniqueInterFI}, \eqref{Ec-Propo53} and \eqref{Ec-s-06}, we can bound $\mu(B(x;r))$ as follows: 
\begin{align}
\mu(B(x;r)) &\leq \mu\left( B(x; |\clC_{n+1}(a_1,\ldots,a_{n+1})|)\right) \nonumber\\
&\leq \mu\left( \bigcup_{j=-1}^{2} \clC_{n+1}(a_1,\ldots,a_{n},a_{n+1} + j)\right) \nonumber\\
&= 4\prod_{k=1}^{n+1} \frac{1}{(N-1)s_k} \leq 4\left( \frac{1}{s_{n+2}(N^{n+1} s_1\cdots s_{n+1})^2}\right)^s. \label{Ec-Caso1-01}
\end{align}
Moreover, our assumption on $r$ and \eqref{Ec-CotaDiamJna} give
\begin{equation}\label{Ec-Caso1-02}
r>|J_{n+1}(\bfa)| \geq \frac{1}{(N^{n+1}s_1\cdots s_{n+1})^2s_{n+2}}. 
\end{equation}
Combining \eqref{Ec-Caso1-01} and \eqref{Ec-Caso1-02}, we conclude $\mu(B(x;r))\leq 4r^s$.

\paragraph{Case II.} Assume that $|\clC_{n+1}(\bfa)|\leq r<|J_n(\bfa)|$. Let us estimate how many fundamental intervals of order $n+1$ intersect $B(x;r)$. If $B(x;r)$ contains $m$ cylinders of order $n+1$, then it contains at most $m+2$ fundamental intervals of order $n+1$ and it intersects at most $m+4$ fundamental intervals of order $n+1$. Hence, using $|B(x;r)|=2r$, $|\clC_{n+1}(\bfa)|\geq N^{-2(n+1)}(s_1\cdots s_ns_{n+1})^{-2}$, and $1\leq rN^{2(n+1)}(s_1\cdots s_n)^2 s_{n+1}$, we have
\[
\frac{m}{N^{2(n+1)}(s_1\cdots s_{n+1})^2} \leq 2r 
\]
and hence
\[
m+4\leq 6rN^{2(n+1)}(s_1\cdots s_{n+1})^2 .
\]

Proposition \ref{PropoUniqueInterFI} implies $\mu(B(x;r))\leq \mu(J_n(\bfa))$; therefore, since $\min\{a,b\}\leq a^sb^{1-s}$ for all $a,b>0$, 
\begin{align*}
\mu(B(x;r)) &\leq \min\left\{ \mu(J_n(\bfa)), 6rN^{2(n+1)}(s_1\cdots s_ns_{n+1})^2  \frac{1}{(N-1)^{n+1}s_1\cdots s_{n+1}}\right\}\nonumber\\
&= \frac{1}{(N-1)^ns_1\cdots s_{n}}\min\left\{ 1,6rN^{2(n+1)}(s_1\cdots s_n)^2s_{n+1}  \frac{1}{(N-1)}\right\}\nonumber\\
&\leq \frac{1}{(N-1)^ns_1\cdots s_{n}} \left( 6r N^{2n} \frac{N}{N-1} (s_1\cdots s_n)^2s_{n+1}\right)^s \nonumber\\
&= \left(\prod_{j=1}^n (N-1)s_j\right)^{-1} \left(s_{n+1} \left( \prod_{k=1}^n Ns_k\right)^2 \right)^s \left( \frac{N}{N-1} \right)^s 6^sr^s, \nonumber
\end{align*}
so, by \eqref{Ec-s-06}, $\mu(B(x;r))< 12r^s$.

We now apply the Mass Distribution Principle to conclude $\dim_H F\geq s_0$.
\subsubsection*{Proof of Proposition \ref{PropoUniqueInterFI}}

Now we prove Proposition \ref{PropoUniqueInterFI}. Keep the notation as above. We consider three cases depending on the value of $a_n$. In each case, we estimate the distance between $J_n(\bfa)$ and the closest fundamental intervals of order $n$.

\paragraph{Case \textsc{I}.} Assume that $s_n+1\leq a_n\leq Ns_n-2$. The fundamental intervals of order $n$ which are neighbors of $J_n(\bfa)$ are $J_n(a_1,\ldots,a_{n-1},a_n-1)$ and $J_n(a_1,\ldots,a_{n-1},a_n+1)$. Put
\[
\bfa':=(a_1,\ldots,a_{n-1},a_n-1), \quad
\bfa'':=(a_1,\ldots,a_{n-1},a_n+1),
\]
so $J_n(\bfa')$ lies to the right of $J_n(\bfa)$ and $J_n(\bfa'')$ to the left. 

The distance between $J_n(\bfa)$ and $J_n(\bfa')$, $d\left(J_n(\bfa'),J_n(\bfa)\right)$, satisfies
\begin{align*}
d\left(J_n(\bfa'),J_n(\bfa)\right) &:= \inf\left\{|y - z|:y\in J_n(\bfa), z\in J_n(\bfa')\right\} \nonumber\\
&=\inf J_n(\bfa') - \sup J_n(\bfa) \nonumber\\
&= S_{\bfa}^{n-1} \circ S_{a_n-1}^1(0) - S_{\bfa}^{n-1}\circ S_{a_{n}}^1 ((s_{n+1}-1)^{-1})) \nonumber\\
&= \left( \prod_{k=1}^{n-1} \frac{1}{a_k(a_k-1)}\right)\left( S_{a_n-1}(0) - S_{a_n}^1 ((s_{n+1}-1)^{-1}) \right) \nonumber\\
&= \left( \prod_{k=1}^{n-1} \frac{1}{a_k(a_k-1)}\right)\left( \frac{1}{a_n-1} - \frac{1}{a_n(a_n-1)(s_{n+1}-1)} - \frac{1}{a_n} \right) \nonumber\\
&= \left( \prod_{k=1}^{n} \frac{1}{a_k(a_k-1)}\right)\left( \frac{s_{n+1}-2}{s_{n+1}-1} \right) \nonumber\\
&= |J_n(\bfa)|(s_{n+1}-2)>|J_n(\bfa)|>r.
\end{align*}
The distance between $J_n(\bfa)$ and $J_n(\bfa'')$ satisfies
\begin{align*}
d(J_n(\bfa''),J_n(\bfa)) &= \inf J_n(\bfa) - \sup J_n(\bfa'') \nonumber\\
& = S_{\bfa}^{n-1} \circ S_{a_n}^1(0) - S_{\bfa}^{n-1}\circ S_{a_{n}+1}^1 \left((s_{n+1}-1)^{-1} \right) \nonumber\\
&= \left( \prod_{k=1}^{n-1} \frac{1}{a_k(a_k-1)}\right)\frac{1}{a_n(a_n+1)} \left( 1 - \frac{1}{s_{n+1}-1} \right) \nonumber\\
&= |J_n(\bfa)| \left(\frac{a_n-1}{a_n+1}\right) (s_{n+1}-2)> |J_n(\bfa)|>r \nonumber.
\end{align*}
The next to last inequality holds because $(a_{n}-1)(a_{n}+1)^{-1}(s_{n+1}-2)>1$ is equivalent to
\[
s_{n+1}>3 + \frac{2}{a_n-1},
\]
which follows from $s_{n+1}\geq 4$ and $a_n\geq 4$.

In a few words, the distance between $J_n(\bfa)$ and its neighboring fundamental intervals of order $n$ is larger than $r$ and the proof for Case I is done.

Note that we can still use the argument involving $\bfa'$ and $\bfa$ when $a_n=Ns_n-1$ and that we can use the argument involving $\bfa''$ and $\bfa$ when $a_n=s_n$.

\paragraph{Case \textsc{II}.} Suppose $a_n=Ns_n-1$. In this case, all the fundamental intervals of order $n$ different from $J_n(\bfa)$ and contained in $\clC_{n-1}(\bfa)$ are to the right of $J_n(\bfa)$. Then, by the last sentence of the previous case, it suffices to estimate $\inf J_n(\bfa) - \inf \clC_{n-1}(\bfa)$:
\begin{align*}
\inf J_n(\bfa) - \inf \clC_{n-1}(\bfa) &= S_{\bfa}^{n-1}\circ S_{a_n}^1(0) - S_{\bfa}^{n-1}(0) \nonumber\\
&= \left(\prod_{k=1}^{n-1} \frac{1}{a_k(a_k-1)}\right) \left(\frac{1}{a_n} - 0\right) \nonumber\\
&> \frac{1}{s_{n+1}-1}\prod_{k=1}^{n} \frac{1}{a_k(a_k-1)} = |J_{n}(\bfa)|>r.
\end{align*}
Hence, the distance between $J_n(\bfa)$ and its neighboring fundamental intervals of order $n$ exceeds $r$.

\paragraph{Case \textsc{III}.} Suppose $a_n=s_n$. We will use the next technical lemma.
\begin{lem01}\label{LemInducMath}
Let $(x_m)_{m\geq 0}, (y_m)_{m\geq 0}$ be two sequences on $ \Na_{\geq 2}$ and define $(z_m)_{m\geq 0}$ by $z_m=x_my_m$. Then, 
\[
\forall m\in\Na_0 \quad \sum_{r=0}^{m} \frac{1}{x_{m-r}}\prod_{j=0}^{m-r} z_j<\prod_{j=0}^{m}z_j.
\]
We also have
\[
\forall m\in \Na_0 \quad \sum_{r=0}^{m} \frac{1}{x_{m-r}}\prod_{j=0}^{m-r} z_j = \sum_{i=0}^m \frac{1}{x_i}\prod_{k=0}^iz_k.
\]
\end{lem01}
\begin{proof}
The first part is shown by induction on $m$ and noting that $z_m=x_my_m> x_m+y_m$ holds because $x_m,y_m\geq 2$. The second part follows from a change of variable on the indexes.
\end{proof}

Let $J_n(\bfb)$, $\bfb=(b_1,\ldots,b_n)$, be the fundamental interval of order $n$ immediately to the right of $J_n(\bfa)$ and let $j\in\{0..n-1\}$ be such that
\begin{align*}
b_1=a_1, &\;\ldots, \; b_{n-j-1} =a_{n-j-1}, \; b_{n-j}= a_{n-j}-1, \nonumber\\
b_{n-j+1} &= N s_{n-j+1} - 1 , \;\ldots, \; b_n=Ns_n-1. \nonumber
\end{align*}
The proof of Proposition \ref{PropoUniqueInterFI} will be complete once we show that
\begin{equation}\label{Caso03-ec01}
\inf J_n(\bfb) - \sup J_n(\bfa)> |J_n(\bfa)|.
\end{equation}
The left-hand side can be rewritten as
\begin{align*}
\inf J_n(\bfb) &- \sup J_n(\bfa) = S_{\bfb}^n(0) - S_{\bfa}^n\left( \frac{1}{s_{n+1}-1}\right) \nonumber\\
&=S_{\bfa}^{n-j-1}\circ S_{b_{n-j}\ldots b_n}^{j+1}(0) - S_{\bfa}^{n-j-1}\circ S_{a_{n-j}\ldots a_n}^{j+1}\left( \frac{1}{s_{n+1}-1} \right) \nonumber\\
&= \left( \prod_{k=1}^{n-j-1}\frac{1}{a_k(a_k-1)}\right) \left(S_{b_{n-j}\ldots b_n}^{j+1}(0) - S_{a_{n-j}\ldots a_n}^j\left( \frac{1}{s_{n+1}-1} \right)\right) \nonumber
\end{align*}
and, in view of \eqref{DiamJna}, the inequality \eqref{Caso03-ec01} follows from
\begin{equation}\label{Caso03-ec02}
S_{b_{n-j}\ldots b_n}^{j+1}(0) - S_{a_{n-j}\ldots a_n}^{j+1}\left( \frac{1}{s_{n+1}-1} \right)
>
\frac{1}{s_{n+1}-1}\prod_{k=n-j}^{n}\frac{1}{a_k(a_k-1)}.
\end{equation}
The terms on the left in \eqref{Caso03-ec02} are
\begin{align*}
S_{b_{n-j}\ldots b_n}^{j+1}(0)  &= \sum_{r=0}^{j-1} (b_{n-r}-1) \prod_{k=r}^{j}\frac{1}{b_{n-k}(b_{n-k}-1)} + \frac{1}{a_{n-j}-1}, \nonumber\\
S_{a_{n-j}\ldots a_n}^{j+1}\left( \frac{1}{s_{n+1}-1} \right)&= \frac{1}{s_{n+1}-1} \prod_{k=0}^{j}\frac{1}{a_{n-k}(a_{n-k}-1)}+ \nonumber\\
&\;\; + \sum_{r=0}^{j-1} (a_{n-r}-1) \prod_{k=r}^{j}\frac{1}{a_{n-k}(a_{n-k}-1)} + \frac{1}{a_{n-j}}; \nonumber
\end{align*}
hence, 
\begin{align*}
S_{b_{n-j}\ldots b_n}^{j+1}&(0) - S_{a_{n-j}\ldots a_n}^{j+1}\left( \frac{1}{s_{n+1}-1} \right) > \frac{1}{a_{n-j}(a_{n-j}-1)} + \nonumber\\
&\quad - \left( \frac{1}{s_{n+1}-1} \prod_{k=0}^{j}\frac{1}{a_{n-k}(a_{n-k}-1)} + \sum_{r=0}^{j-1} (a_{n-r}-1) \prod_{k=r}^{j}\frac{1}{a_{n-k}(a_{n-k}-1)}\right). \nonumber
\end{align*}
Inequality \eqref{Caso03-ec02} will be proven if we show that
\[
\frac{1}{a_{n-j}(a_{n-j}-1)} -\sum_{r=0}^{j-1} (a_{n-r}-1) \prod_{k=r}^{j}\frac{1 }{a_{n-k}(a_{n-k}-1)}> \frac{2}{s_{n+1}-1} \prod_{k=0}^{j}\frac{1}{a_{n-k}(a_{n-k}-1)},
\]
which is
\[
\prod_{k=0}^{j-1} a_{n-k}(a_{n-k}-1) - \sum_{r=0}^{j-1} \frac{1}{a_{n-r}} \prod_{k=0}^{r} a_{n-k}(a_{n-k}-1)> \frac{2}{s_{n+1}-1}.
\]
Finally, considering $m=j-1$, $x_k=a_{n-k}$, $y_k=a_{n-k}-1$, and $z_k=x_ky_k$ for $k\in \{0,\ldots,m\}$, the previous inequality becomes
\[
\prod_{k=0}^{m} z_k - \sum_{r=0}^{m} \frac{1}{x_{r}} \prod_{k=0}^{r} z_{k}> \frac{2}{s_{n+1}-1},
\]
which holds by Lemma \ref{LemInducMath} and $s_{n+1}-1\geq 3$. Therefore, \eqref{Caso03-ec02} follows and the proof of Proposition \ref{PropoUniqueInterFI} is complete.

\section{Final Remarks}

While our strategy for the lower bound on Theorem \ref{TEOMAIN} resembles that in \cite{suwu14}, our argument for the upper bound does not. We had to adopt a new strategy because, in general, Lüroth series do not converge as fast as regular continued fractions. To be more precise, if $x=\left\langle a_1,a_2,\ldots\right\rangle \in(0,1]$, $(P_n)_{n\geq 1}$, $(Q_n)_{n\geq 1}$ are as before and $(p_n)_{n\geq 0}$, $(q_n)_{n\geq 0}$ are the numerators and denominators of the regular continued fractions convergents of $x$, then
\[
\forall n\in\Na \qquad | Q_nx - P_n|<\frac{1}{a_n-1},\; |q_nx-p_n|<\frac{1}{q_n}.
\]
In general, while the sequence $(|q_nx-p_n|)_{n\geq 0}$ converges to $0$, $(|Q_nx-P_n|)_{n\geq 1}$ might not (in our context, however, $(|Q_nx-P_n|)_{n\geq 1}$ does converge to $0$).

Our proof of the upper bound in \eqref{Eq-MainResult} also gives half of Theorem \ref{TeoSunWu}. We sketch the argument and leave the details to the reader. Given a real number $x$, denote its regular continued fraction by $[a_0;a_1,a_2,\ldots]$. 

First, for any $\bfb=\sabu\in \Na^{\Na}$ and $n\in\Na$ write
\[
\clC^{cf}_n(\bfb):=\{ x=[0;a_1,a_2,\ldots]\in[0,1]\setminus\QU: a_1=b_1,\ldots,a_n=b_n\}.
\]
Fix $\alpha>0$ and define
\[
G^{cf}(\alpha) = \left\{x=[0;a_1,a_2,\ldots]\in [0,1]\setminus\QU: \lim_{n\to\infty} \frac{\log q_n}{\log a_{n+1}}= \alpha\right\}.
\]
Let $\veps>0$ be an arbitrary positive number and put
\[
C_{\veps}(\alpha):= \left\{ x=[0;a_1,a_2,\ldots]: \limsup_{n\to\infty} \frac{\log q_n}{\log a_{n+1}} \leq \alpha + \frac{\veps}{2}\right\},
\]
so $G^{cf}(\alpha)\subseteq C_{\veps}(\alpha)$. Adapt to the continued fraction context the definitions of $C_{\veps}^n(\alpha)$ and $K_{n}(\bfb)$ for $n\in\Na$ and $\bfb\in\Na^n$. Take $\delta>0$ small and set $s=(2+ (\alpha+\veps)^{-1})^{-1} + \delta$. Using the theory of regular continued fractions (v.gr. \cite{khin}), we can show that for every large $n\in\Na$, we have that
\[
\forall\bfa\in\Na^n \qquad \sum_{b} |K_n(\bfa b)|^s\leq |K_n(\bfa)|^s,
\]
where the sum runs along the set $\{b\in \Na: b^{\alpha+\veps}\geq q_n\}$. As a consequence, by Lemma \ref{LeGJL02}, we have
\[
\dim_H G^{cf}(\alpha)\leq  \frac{1}{ 2 + \frac{1}{\alpha+\veps}} + \delta,
\]
and, since $\delta>0$ can be arbitrarily small, $\dim_H G^{cf}(\alpha)\leq \left( 2 + \tfrac{1}{\alpha+\veps}\right)^{-1}$. Finally, since $x\mapsto (2+ x^{-1})^{-1}$ is increasing and continuous on $\RE_{>0}$ and $\veps>0$ was arbitrary, we conclude that 
\[
\dim_H G^{cf}(\alpha) \leq \frac{1}{ 2 + \frac{1}{\alpha}}.
\]
Writing $\beta=\alpha^{-1}$, the previous inequality becomes
\[
\dim_H \left\{ x=[0;a_1,a_2,\ldots]\in (0,1): \lim_{n\to\infty} \frac{\log a_{n+1}}{\log q_n} = \beta\right\} \leq \frac{1}{2+\beta},
\]
which is half of Theorem \ref{TeoSunWu}.

\Addresses

\begin{thebibliography}{99}
\bibitem {bipe} Bishop, C.; Peres, Y., \textit{Fractals in probability and analysis.}
Cambridge Studies in Advanced Mathematics, 162. Cambridge University Press, Cambridge, 2017. 

\bibitem {DaKra02} Dajani, K.; Kraaikamp, C., \textit{Ergodic theory of numbers.} Carus Mathematical Monographs, 29. Mathematical Association of America, Washington, DC, 2002.

\bibitem {DaKra} Dajani, K.; Kraaikamp, C., \textit{On approximation by Lüroth series.} J. Théor. Nombres Bordeaux 8 (1996), no. 2, 331--346.

\bibitem {faliwawu} Fan, A.-H.; Liao, L.-M; Wang, B.-W.; Wu, J., \textit{On Khintchine exponents and Lyapunov exponents of continued fractions.} Ergodic Theory Dynam. Systems 29 (2009), no. 1, 73--109.

\bibitem {ftz} Feng, Y; Tan, B.; Zhou, Q.-L., \textit{Exact Dimensions of Exceptional Sets in Lüroth Expansions}, Fractals, doi: 10.1142/S0218348X21501425


\bibitem {gero-good} Gonzalez Robert, G., \textit{Good's theorem for Hurwitz continued fractions}.  Int. Journal of Number Theory. (2020) Volume No.16, Issue No. 07.

\bibitem {good} Good, I. J., \textit{The fractional dimensional theory of continued fractions.} Proc. Cambridge Philos. Soc. 37 (1941), 199--228.

\bibitem {jagervroedt} Jager, H.; de Vroedt, C., \textit{Lüroth series and their ergodic properties.} Nederl. Akad. Wetensch. Proc. Ser. A 72=Indag. Math. 31 1969 31--42.

\bibitem {khin}  Khinchin, A. Ya., \textit{Continued fractions.} Translated from the third (1961) Russian edition. Reprint of the 1964 translation. Dover Publications, Inc., Mineola, NY, 1997.

\bibitem {luroth} Lüroth, J., \textit{Ueber eine eindeutige Entwickelung von Zahlen in eine unendliche Reihe.} Math. Ann. 21 (1883) 411--423.

\bibitem {mun11} Munday, S., \textit{A note on Diophantine fractals for $\alpha$-L\"uroth systems.} Integers 11B (2011), Paper No. A10, 14 pp.

\bibitem {parth} Parthasarathy, K. R. \textit{Introduction to probability and measure.} Texts and Readings in Mathematics, 33. Hindustan Book Agency, New Delhi, 2005. 

\bibitem {suwu14} Sun, Y.; Wu, J. \textit{A dimensional result in continued fractions.} Int. J. Number Theory 10 (2014), no. 4, 849--857.

\end{thebibliography}
\end{document}